\documentclass[english]{amsart}
\usepackage{amscd,amsfonts,amssymb,amsmath,a4wide}

\def\proof {\noindent{\sc{ Proof. }}}
\def\QED {\mbox{}\hfill {\small \fbox{}} \\}  

\def\Zm{{\mathbb Z}}

\def\Rm{{\mathbb R}}
\def\Tm{{\mathbb T}}
\def\Nm{{\mathbb N}}

\def\lto{\longrightarrow}
\def\lmto{\longmapsto}
\def\eq{\Longleftrightarrow}
\def\leq{\leqslant}
\def\geq{\geqslant}

\newcommand{\mE}{\ensuremath{\mathcal{E}}}

\newcommand{\vs}{\vspace{.2cm}}
\makeatletter
\renewcommand{\section}{\@startsection
{section}
{1}
{0mm}
{-1.2\baselineskip}
{\baselineskip}
{\center \scshape}}

\renewcommand{\subsection}{\@startsection
{subsection}
{2}
{0mm}
{-\baselineskip}
{0mm}
{\normalfont \normalsize \bfseries}}
\makeatother

\author{Patrick Bernard}
\address{Institut Fourier
BP 74,
 38402, Saint Martin d'H\`eres Cedex,
France}
\email{
Patrick.Bernard@ujf-grenoble.fr}
\urladdr{http://www-fourier.ujf-grenoble.fr/~pbernard/}
\title{The asymptotic behavior
of solutions of  forced Burgers equation on the
circle}
\begin{document}

\begin{abstract}
We describe the asymptotic behaviour of entropy solutions
of unviscid Burgers equation on the circle
with time-periodic forcing.
These solutions converge to periodic states, but the period
of these limit states may be greater than the period of the forcing.
We obtain as a corollary a new result on the flow of the associated
Hamiltonian system.
\end{abstract}

\maketitle

\section{Introduction}
\subsection{}
The standard circle $\Rm/\Zm$ is noted $\Tm$.
The cotangent bundle $T^*\Tm$ is identified
with $\Tm\times \Rm$.
Given a function $f(t,x)$ of two variables,
we will note
$f_t$ the function
$x\lmto f(t,x)$.
The partial derivative with respect
to the variable $t$ will be denoted $\partial_t f$.
In all this paper, we will consider a 
time-periodic Hamiltonian 
$H(t,x,p):\Rm\times T^* \Tm
=\Rm\times \Tm \times \Rm \lto \Rm$,
the associated time-periodic vector-field
of 
$
\Tm \times \Rm
$
is noted $X$.
We have 
$$
X(t,x,p)=\big(
\partial_p H(t,x,p),
-\partial_x H(t,x,p)\big).
$$

\subsection{}\label{Hhyp}
The following standard hypotheses will be assumed :\\
$i.$ The Hamiltonian $H$ is $C^2$ and $1$-periodic in $t$.\\
$ii.$ The Hamiltonian $H$ is convex in $p$, and
$\partial_{pp}H>0$.\\
$iii.$ The Hamiltonian has superlinear growth in $p$, 
\textit{i.e.}
$\lim_{|p| \lto \infty}H(t,x,p)/p=\infty$ 
for each $(t,x)$.\\
$iv.$ The Hamiltonian flow is complete. 
More precisely
for all $(t_0,x_0,p_0)$, there exists a $C^1$ curve
$\gamma(t)=(x(t),p(t)):\Rm\lto \Tm\times \Rm$
such that $(x(t_0),p(t_0))=(x_0,p_0)$
and $\dot \gamma(t)=X(t,\gamma(t))$ for all $t\in\Rm$.
The mapping $\gamma(t_0)\lmto \gamma(t)$
is a diffeomorphism of $\Tm\times\Rm$, denoted
$\phi_{t_0,t}$. We will pay a special attention 
to the diffeomorphism $\phi=\phi_{0,1}$.
Note that the completeness Hypothesis is satisfied if
there exists a constant $C$ such that
$
|H_t|\leq C(1+H)$.

\subsection{}
A typical example of Hamiltonian satisfying our
hypotheses is
$$
H(t,x,p)=\frac{1}{2}p^2 +V(t,x)
$$
with a $C^2$ potential $V$ periodic in $t$.

\subsection{}\label{burgers}
We consider the equation
$$
\hspace{5.8cm}
\partial_t y+\partial_x \big( H(t,x,y) \big)=0
\hspace{5cm}(B)
$$
of the unknown function 
$y(t,x): \Rm\times \Tm \lto \Rm
$.
This equation will be called the 
Burgers equation in the sequel.
Note that in  case $H=p^2/2+V(t,x)$,
we have the standard 
forced unviscid Burgers equation
$$
\partial_t y+y\partial_x y =-\partial_x V(t,x).
$$

\subsection{}\label{char}
The Burgers equation is quasi-linear,
and its characteristics 
(see \cite{Arnold}, chapter 2)
are the trajectories of $X$.
In other words,
if $y(t,x):[a,b]\times \Tm\lto \Rm$
is a $C^1$ solution of the Burgers
equation,
then  for each $a\leq t_0\leq t_1\leq b$,
the graph of the function $x\lmto y(t_1,x)$
is the image by the 
diffeomorphism  $\phi_{t_0,t_1}$ of  
the graph of the function 
$x\lmto y(t_0,x)$.

\subsection{}\label{HJ}
Still assuming that $y$ is a  $C^1$ solution
of the Burgers equation,
we obtain that
$$
c(t)=\int_{\Tm} y(t,x)dx
$$
is a constant, that we denote $c$.
The function $y$ can be written
$y(t,x)=c+\partial_x u(t,x),$
where $u(t,x):I\times \Tm\lto \Rm$
satisfies the Hamilton-Jacobi equation
$$
\hspace{5.5cm}
\partial_t u+H(t,x,c+\partial_x u) =0.
\hspace{4.5cm}
(HJc)
$$

\subsection{}
It is known that there exists in general no 
classical solution
of Burgers equation
defined
on $\Rm^+\times \Tm$
satisfying a given initial condition
$y(0,x)=y_0(x)$.
However, the Cauchy problem is well-posed
in the sense of entropy solutions.
More precisely,
for all $s\in\Rm$ and all function
$u_s\in L^{\infty}(\Tm)$ there exists
a unique entropy solution
$u(t,x)\in C([s,\infty),L^1(\Tm))$
such that $u(s,.)=u_s$, see \ref{propentropy}.
The operator which, to each function 
$u_0$, associates the function 
$u_1=u(1,.)$, where $u(t,x)$
is the entropy solution with initial
condition $u(0,.)=u_0$,
can be extended to a continuous operator 
$E:L^1(\Tm)\lto L^1(\Tm)$, see \ref{Eprop}.

\subsection{}
We want to describe the asymptotic behaviour
of entropy  solutions.
Let us first recall that for each $c$ there exists
a $1$-periodic solution of average $c$.
More precisely, there exists a continuous
and increasing function $c\lmto y^c$ from $\Rm$ to $L^1(\Tm)$
which, to each $c$, associates
a fixed point of $E$ of average $c$, see 
\cite{KO} or \cite{JKM}.
Note however that there may exist more than one
fixed point of a given average $c$.
It is natural to ask whether all solutions are 
attracted by these fixed points.
The answer is negative, 
 there are examples
where there exist periodic points of $E$ 
which are not fixed point, see \cite{FM},
that is periodic entropy 
solutions of minimal period greater than one.
These subharmonic solutions in turn attract
all other solutions, as we now state.

\subsection{}\label{main}
\textsc{Theorem}
\begin{itshape}
Let $y(t,x): [t_0,\infty)\times \Tm\lto \Rm$
be an entropy solution
of Burgers equation.
There exists an integer $T$
and 
an entropy solution $\omega(t,x):\Rm\times \Tm\lto \Rm$
which is $T$-periodic in $t$ and such that
$$
y_{t+nT}\lto \omega_t
$$
in $L^1(\Tm)$
as $n$ goes to infinity.
\end{itshape}
\vs\\
If $H$ is a function of $p$ only,
then $\omega(t,x)$ is the constant
$\int y(t,x)dx$. The result in this special
case has been obtained by Lax \cite{Lax}.
If $H$ does not depend on $t$, then the asymptotic
solution $\omega$ does not depend on $t$ either,
the result in this case follows from 
works of Roquejoffre \cite{Roque} and Fathi \cite{F3}.
The theorem will be proved in section \ref{E}
as a consequence of a similar result for viscosity solutions
of $(HJc)$
obtained in \cite{Be}, see also \cite{BR}.

\subsection{}
One can compare the situation with 
the viscous case. 
If one considers the parabolic equation
$$
\hspace{5.5cm}
\partial_t y+\partial_x \big( H(t,x,y) \big)=\mu \partial_{xx}y
\hspace{4.5cm}(B_{\mu})
$$
with $\mu>0$, the behaviour is much simpler.
One can prove in the line of \cite{JKM}, see also \cite{Bessi},
that for each $c$, there exists a unique solution
$y^c$ of average $c$ which is $1$-periodic in time.
This solution attracts all the solutions of average $c$.
More precisely, if $y:[t_0,\infty)\times \Tm \lto \Rm$
is a solution of $(B_{\mu})$, and if 
$\int y_t dx =c$, then
$y_{t+n}\lto y^c_t $ uniformly
as $n\in \Nm \lto \infty$.

\subsection{}
The result in the unviscid case  can be used to study
the dynamics of the diffeomorphism 
$\phi=\phi_{0,1}$.
Note that this diffeomorphism is
a finite composition of area preserving 
right twist maps, and that 
any finite composition of
right twist maps can be obtained that way.
This  correspondence between the twist property
and the convexity of the Hamiltonian
has been described by Moser, see \cite{Moser}.

\subsection{}\label{graph}
In order to give a more geometrical meaning to Theorem \ref{main},
we consider the set $\mE$ of functions
$f:\Tm\lto \Rm$ which can be locally 
written as the sum of a continuous and of a decreasing function.
A function $f\in \mE$
has  a right limit 
$f^-(x)$ and a left limit
$f^+(x)$ at each point $x$.
These limits  satisfy 
$f^-(x)\leq
f^+(x)$,
with a strict inequality on an at most countable
set.
Let $G^-(f)$ and $G^+(f)$ be the graphs,
in $\Tm\times\Rm$ of the functions $f^-$ and
$f^+$.
We define the graph  $G(f)=G^-(f)\cap G^+(f)$
of $f$.
Note that $\overline {G(f)}=G^-(f)\cup G^+(f)$.
It is also useful to consider
the set 
$$
H(f)=\bigcup_x \{x\}\times[f^-(x),f^+(x)]\subset \Tm\times \Rm,
$$
which is a Jordan curve containing $\overline {G(f)}$.
The Hausdorff distance  $d_H(f,g)$
between the compact sets $H(f)$ and $H(g)$
defines a distance  $d_H$ on $\mE$
(one should take the quotient of $\mE$ by the relation
of almost everywhere equality).

\subsection{}
The link between entropy solutions and the dynamic 
of $\phi$ can now be detailed, see \ref{Edyn}.
For each $y\in L^1$, we have $E(y)\in \mE$.
If in addition $y\in \mE$, then 
$$
\overline{G(E(y))}\subset \phi(G(y)).
$$
This property has striking consequences.
For example, if $y^c$ is a fixed point of $E$,
then $y^c\in\mE$ and 
$\overline{G(y^c)}$ is negatively invariant.
As a consequence, the $\alpha$-limit of $\phi_{|G(y^c)}$
is a non-empty compact set which is fully invariant by $\phi$.
It is an Aubry-Mather set.
The rotation number $\rho(c)\in \Rm$ of the orbits
of this set depends only on $c$, and the function
$c\lmto \rho(c)$ is non decreasing and continuous,
see \ref{Erotation}.

\subsection{}\label{main2}
Having defined the rotation number $\rho(c)$
allows us to complement Theorem \ref{main}.
The asymptotic behaviour of solutions
depend strongly on their space average $c$ and on the associated
rotation number $\rho(c)$.
If $\rho(c)$ is irrational, then there exists a single fixed
point of $E$ of average $c$, see \cite{WE}.
We will prove that it attracts all the
trajectories of average $c$, that is one can take $T=1$ in Theorem
\ref{main}.
If $\rho(c)$ is rational, 
$\rho(c)=p/q$ in lowest terms,
then one can take $T=q$ in Theorem \ref{main}.
It is thus natural to define
the  integer $T(c)$ by $T(c)=1$
if $\rho(c)$ is irrational and $T(c)=q$ if 
$\rho(c)=p/q$ in lowest terms,
and we have the following refinement of Theorem \ref{main} :\vs\\
\textsc{Theorem}
\begin{itshape}
For each $y\in L^1(\Tm)$,
and $c=\int y$,
there exists a fixed point $l\in \mE$
of $E^{T(c)}$
such that 
$\|E^{nT(c)}(y)-l\|_{L^1}\lto 0$, 
which implies that 
$d_H(E^{nT(c)}(y),l)\lto 0$. 
\end{itshape}
\subsection{}\label{cor}
We obtain a new result on the dynamics of $\phi$,
which may be seen as a converse to the 
celebrated result of Birkhoff
(see \cite{Mathertwist},  \cite{HF}, \cite{Siburg})
stating that a rotational 
(not homotopic to a constant) 
Jordan curve in $\Tm\times \Rm$ 
which is invariant by $\phi$ 
has to be the graph of a Lipshitz function
$y:\Tm\lto \Rm$.
\vs\\
\textsc{Corollary}
\begin{itshape}
Let $G\subset \Tm\times\Rm$ be the graph 
of a continuous function $y:\Tm\lto\Rm$.
Assume that there exists an increasing sequence
$n_k$ of positive integers such that 
$\phi^{n_k}(G)$ is the graph of a continuous function.
Then there exists a positive integer $T$
such that $\phi^T(G)=G$,
and the function $y$ is Lipschitz.
In addition, if $\phi(G)\neq G$,
then the rotation number of $\phi^T_{|G}$ is an integer,
hence $G$ contains a $T$-periodic point of $\phi$.
\end{itshape}
\vs\\
Let us mention that it is certainly possible to give a more direct proof
of this Corollary.
One could use a topological approach, as suggested to me by 
by P. Le Calvez 
or a variational approach, as suggested by J. Xia.
The proof presented here in \ref{reversed} and \ref{coroproof}
is however extremely short.

\subsection{}
In the rest of the paper, we will detail the outline 
given above.
We will obtain all the important properties of entropy  solutions
of $(B)$
as consequences of properties of 
the viscosity solutions of $(HJc)$.
Hence we first describe these viscosity solutions
in section \ref{V}, and draw our conclusions in section \ref{E}.

\section{Calculus of variations and Hamilton-Jacobi equation}\label{V}
In  the present section, we describe the
main properties of viscosity solutions of the equation $(HJc)$.
These properties follow from the study of extremals
via the Hopf-Lax-Oleinik formula, a global reference is the work 
of Fathi, \cite{F1} and  \cite{Fbook}.
We also state a result analogous to Theorem \ref{main}
for these solutions.

\subsection{}\label{Lhyp}
It is  useful to introduce the Lagrangian
$L:\Rm\times \Tm\times \Rm\lto \Rm$
associated to $H$.
It is defined by
$$
L(t,x,v)=\sup _p pv-H(t,x,p),
$$
and has the following properties,
which follow easily from the 
analogous properties \ref{Hhyp}
 of $H$:\\
$i.$ The Lagrangian  $L$ is $C^2$ and $1$-periodic in $t$.\\
$ii.$ The Lagrangian  $L$ is convex in $v$, and
$\partial_{vv}L>0$.\\
$iii.$ The Lagrangian  has superlinear growth in $v$, 
\textit{i.e.}
$\lim_{|v| \lto \infty}L(t,x,v)/v=\infty$ 
for each $(t,x)$.\vs\\
The Lagrangian associated to the modified Hamiltonian
$H(t,x,p+c)$ is $L(t,x,v)-cv$,
it satisfies the three properties above.

\subsection{}
For each $c$ and each
$t_0\leq t$,
we have the Hopf-Lax-Oleinik
operator
$V^c _{t_0,t}:C(\Tm,\Rm)\lto C(\Tm,\Rm)$
defined by
$$
V^c_{t_0,t}(u)(x)=
\min 
\left(
u(x(t_0))
+\int _{t_0}^{t}
L(s,x(s),\dot x(s)) -c\dot x(s)\,ds
\right) 
$$
where the minimum is taken 
on the set of absolutely continuous 
curves $x:[t_0,t]\lto \Tm$ such that
$x(t)=x$.
Any curve realising the above minimum
is $C^2$ and is the projection of a trajectory of $X$.
More precisely, the curve
$(x(s),\partial_v L(s,x(s),\dot x(s)))$
is a trajectory of $X$.
These operators clearly satisfy the Markov property
$$
V^c_{t_1,t_2}\circ
V^c_{t_0,t_1}=
V^c_{t_0,t_2}.
$$
Note that this operator has been used in the study
of viscosity solutions for quite a long time, see for instance
\cite{Fleming}.

%let us define the function 
%$F^c _{t,t'}:\Tm\times \Tm\lto \Rm$
%by 
%$%$
%F^c _{t,t'}(x,x')=
%\min 
%\left(
%\int _t ^{t'}
%L(t,x(s),\dot x(s)) -c\dot x(s)\,ds
%\right) 
%$$
%where the minimum is taken on the set of 
%absolutely continuous curves $x(s):[t,t']\lto \Tm$
%such that $x(t)=x$ and $x(t')=x'$.

\subsection{}\label{Vprop}
A function is called semi-concave if it 
is locally the sum of a 
smooth function and of a concave function.
A function $u(x):\Tm\lto \Rm$
is called $K$-semi-concave
if 
$\partial_{xx}u\leq K$
in the sense of distributions.
\vs\\
\textsc{Proposition}
\begin{itshape}
The following properties are equivalent for a  function 
$u\in C([t_0,t_1]\times \Tm, \Rm)$ : \\
i. 
The function $u$ is a viscosity solution of $(HJc)$
in the classical sense, (see \cite{Barles}).\\
ii. 
The function $u$ is locally Lipschitz on $]t_0,t_1[\times\Tm$ and it
satisfies $(HJc)$ almost everywhere.
In addition, there exists a non increasing function
$K:]t_0,t_1]\lto ]0,\infty[$ such that 
the functions $u_t$ is $K(t)$-semi-concave for $t_0<t \leq t_1$.\\
iii.
For each $t$ and $t'$ such that $t_0\leq t\leq t'\leq t_1$,
we have $u_{t'}=V^c_{t,t'}(u_t)$.
\end{itshape}
\vs\\
\textsc{Corollary}
\begin{itshape}
For each initial condition $u_{t_0}$, 
there exists one 
and only one viscosity solution 
$u:[t_0,\infty) \times \Tm \lto \Rm$
of $(HJc)$,
it is given by 
$$
u(t,x)=
 V^c_{t_0,t}u_{t_0}(x).
$$
There exists a non-increasing function
$K(t):]t_0,\infty[\lto]0,\infty[$
such that  $u_t$ is $K(t)$-semi-concave 
and such that $u$ is $K(t)$-Lipschitz on $[t,\infty)\times \Tm$
for each $t>t_0$.
In addition, if $u_{t_0}$ is Lipschitz, then so is $u$
on $[t_0,\infty)\times\Tm$.
\end{itshape}
\vs\\
\proof
It is standard that $i.\eq iii.$,
see a good exposition in \cite{FMa}. Let us recall a possible  sketch
of  proof.
One can first prove using variations around the
maximum principle that there is at most one function 
satisfying $i.$ with a given initial condition $u_{t_0}$
(see \cite{Barles}, 2.4.).
On the other hand, it is
obvious that there exists  one and only
one function satisfying $iii.$, namely 
$(t,x)\lmto V^c_{t_0,t}u_{t_0}(x)$.
One can prove 
(see \cite{Fbook})
that this function also satisfies $i$.
It is then the only one to do so, by uniqueness.
It is also classical that $iii.\Longrightarrow ii.$,
see  \cite{Fbook}.
We shall prove more carefully that $ii.\Longrightarrow i.$,
which seems less classical.
Let us fix $(S,Q)\in ]t_0,\infty)\times \Tm$,
it is enough to prove
(see \cite{Barles}, 5.3) that  all $C^1 $ function $\phi$
such that $u-\phi$ has a local minimum at $(S,Q)$,
satisfies the equation at $(S,Q)$.
If such a function $\phi$ exists, then 
$\partial_x u (S,Q)$ exists and is equal to $\partial_x \phi(t,x)$.
It follows from the Lemma below that 
$u$ is differentiable at $(S,Q)$, and satisfies
the equation at this point,
which implies that $\phi$ also satisfies the equation
at $(S,Q)$.
The additional conclusions of the Corollary
follow in a classical way from the analysis of calibrated curves,
as defined in \ref{cal},
see \cite{Fbook}. 
\QED
\vs\\
\textsc{Lemma}
\begin{itshape}
Let $u(t,x)$ be a function satisfying condition $ii.$
of the proposition.
If $(S,Q)\in ]t_0,\infty)\times\Tm$
is a point where 
$\partial _x u$ exits, then  the function $u$ is differentiable
and satisfies $(HJc)$ at $(S,Q)$.
\end{itshape}
\vs\\
\proof
Let us fix a time $t_2\in ]t_0,S[$.
In view of the fact that all the functions $u_t, t\geq t_2$
are $K(t_2)$-semi-concave,
it is not hard to prove that 
$\partial_x u(s_n,q_n)\lto \partial_x u(T,Q)$
when $(s_n,q_n)$ is a sequence of points of differentiability of $u$
converging to $(S,Q)$.
If we assume in addition that $(HJc)$ holds at $(s_n,q_n)$,
we obtain that $\partial _t u (s_n,q_n)$
has a limit $H(t,x,c+\partial_x u(t,x))$.
Let us denote $L$ the linear form
$(s,q)\lmto q\partial_x u(S,Q)+sH(S,Q,\partial_x u(S,Q))$.
We have proved that 
there exists a modulus of continuity 
$\delta$ and a set $K\subset \Rm\times \Tm$
of full measure in a neighbourhood
of $(S,Q)$ such that,
for each $(S+s,Q+q)\in K$, the function $u$ is differentiable at
$(S+s,Q+q)$ and 
$\|du(S+s,Q+q)-L\|\leq \delta\big(\|(s,q)\|\big)$.
It follows that we have the estimate
$
|u(S+s,Q+q)-u(S,Q)-L(s,q)|
\leq \|(s,q)\|\delta( \|(s,q)\|)$
for all $(s,q)$ small enough, hence
$u$ is differentiable at $(S,Q)$,
and $du(S,Q)=L$.
\QED

\subsection{}\label{cal}
Let  $u(t,x):[t_0,\infty)\times \Tm\lto \Rm$
be a viscosity solution of $(HJc)$,
and let $t_0\leq t<t'$.
An absolutely continuous  curve 
$x(s):[t,t']\lto \Tm$ is said calibrated by $u$
if 
$$
u(t',x(t'))=
u(t,x(t))+
\int_t^{t'}
L(s,x(s),\dot x(s))-c\dot x(s)\,ds.
$$
If $x(s)$ is a calibrated curve,
then it is $C^2$ and the curve
$\big(x(s),\partial_vL(s,x(s), \dot x(s)\big)$
is a trajectory of the Hamiltonian
vector field $X$.
By extension, we say that a curve
$\gamma(s)=(x(s),p(s)):[t,t']\lto \Tm\times\Rm$
is calibrated by $u$ if
$x(s)$ is calibrated by $u$ and if 
$p(s)=\partial_v L(s,x(s),\dot x(s))$.
It is then a trajectory of $X$.
A curve is said to be calibrated by $u$
on an interval $I$  it is calibrated by $u$
on $[t,t']$
for all $[t,t']\subset I$.
It was proved by Fathi that
if $\gamma(s)=(x(s),p(s)):[t,t']\lto \Tm\times\Rm$
is calibrated, then the function
$u$ is differentiable at $(s,x(s))$ for $s\in ]t,t'[$
and satisfies 
$$
\partial_x u(t,x(s))=p(s).
$$
It is not hard to prove in the same way that
if $u_0$ is Lipschitz, then there exists
a constant $K$ such that 
each calibrated curve 
$\gamma(s)=(x(s),p(s)):[t_0,t]\lto \Tm\times\Rm$
with $t>t_0$ satisfies $|p(s)|\leq K$
for $s\in [t_0,t]$.
We obtain from this remark that $u$ is Lipschitz
if $u_0$ is.

\subsection{}
We will note 
$V^c:C(\Tm,\Rm)\lto C(\Tm,\Rm)$
the operator $ V^c_{0,1}$.
\vs\\
\textsc{Proposition}
\begin{itshape}
The operator $V^c$ is a contraction,
$$\|V^c(u)-V^c(v)\|_{\infty}
\leq
\|u-v\|_{\infty}.$$
In addition,
For each $C>0$, there is a constant
$K>0$ such that, if $|c|\leq C$,
the elements of 
the image of $V^c$
are 
$K$-Lipschitz and $K$-semi-concave.
\end{itshape}
\vs\\
It will also be useful to consider 
the operator $\tilde V:\Rm\times C(\Tm)\lto\Rm\times C(\Tm)$
defined by
$$
\tilde V(c,u)=(c,V^c(u)).
$$
This operator is continuous and compact.

\subsection{}\label{alpha}
There exists a $C^1$ convex and super linear function
$\alpha(c):\Rm\lto\Rm$ such that,
for all solution 
$u(t,x):[t_0,\infty)\times \Tm\lto \Rm$
of $(HJc)$, the function $u(t,x)+t\alpha(c)$
is bounded on $[t_0,\infty)\times \Tm$.
This function is the alpha function 
of Mather,  also called the effective Hamiltonian.
Let us illustrate a bit more the meaning of this function:
\vs\\
%%%%%
%
%
%
%
%
%
%
\textsc{Lemma }
\begin{itshape}
If $u(t,x):[t_0,\infty)\times \Tm\lto \Rm$
is a viscosity solution of $(HJc)$,
the function
$v(\theta,x):\Tm\times \Tm \lto \Rm$
defined by
$$
v(\theta ,x)=\liminf_{t \text{ mod } 1=\theta }
\big(u(t,x)-t\alpha(c)
\big)
$$
is a viscosity solution of the Hamilton-Jacobi 
equation
$$
\partial_{\theta} v
+
H(\theta,x,c+\partial_x v)=\alpha(c),
$$ 
where we have noted $\theta$ for $t \text{ mod }1$.
\end{itshape}
\vs\\
\textsc{Corollary}
\begin{itshape}
The Hamilton-Jacobi equation
$$
\partial_{\theta} u
+
H(\theta,x,c+\partial_x u)=a
$$ 
has a $1$-periodic viscosity solution
if and only if $a=\alpha(c).$
\end{itshape}
\vs\\
A proof of the Lemma is written in \cite{BR}, the corollary
is by now classical.

\subsection{}\label{rotationu}
The number $\rho(c)=\alpha'(c)$
has an important dynamical meaning:\\
\textsc{Proposition }
\begin{itshape}
For each viscosity solution 
  $u(t,x):[t_0,\infty)\times \Tm\lto \Rm$ of $(HJc)$,
there exist curves $\gamma:[t_0,\infty)\lto \Tm$
which are calibrated by $u$.
These curves all have the same rotation number
$$
\lim \frac{1}{t}\int _{t_0}^t \dot x=\rho(c)=\alpha'(c).
$$
\end{itshape}

\proof
Let $u_c:[t_0,\infty)\times \Tm\lto \Rm$ be a 
viscosity solution of $(HJc)$
and let $x(t):[t_0,\infty)\lto \Tm$,
be calibrated by $u_c$.
We have 
$$
u_c(t,x(t))-u_c(t_0,x(t_0))=
\int_{t_0}^t L(s,x(s),\dot x(s)) ds
-c\int_{t_0}^t \dot x(s) ds.
$$
In view of the definition of $\alpha(c)$ in \ref{alpha},
we obtain
$$
\lim _{t\lto\infty }\Big(
\frac{c}{t}\int_{t_0}^t\dot x(s) ds-
\frac{1}{t}
\int_{t_0}^t L(s,x(s),\dot x(s)) ds 
\Big)=\alpha(c).
$$
Let us now consider a viscosity solution
$u_e:[t_0,\infty)\times \Tm\lto \Rm$  of $(HJe)$.
With the same curve $x(s)$, we have 
$$
u_e(t,x(t))-u_e(t_0,x(t_0))\leq
\int_{t_0}^t L(s,x(s),\dot x(s)) ds
-e\int_{t_0}^t \dot x(s) ds
$$
hence 
$$
\alpha(e)
\geq \limsup
\Big(
\frac{e}{t}\int_{t_0}^t\dot x(s) ds
-\frac{1}{t}
\int_{t_0}^t L(s,x(s),\dot x(s)) ds 
\Big),
$$
which implies 
$$
\alpha(e)-\alpha(c)\geq
\limsup \Big(
\frac{e-c}{t}\int_{t_0}^t\dot x(s) ds
\Big).
$$
We conclude
$$
\alpha'(c)\leq \liminf \frac{1}{t}\int_{t_0}^t\dot x(s) ds
\leq  \limsup  \frac{1}{t}\int_{t_0}^t\dot x(s) ds 
\leq \alpha'(c).
$$
It follows that 
the curve  $x(s)$ has a rotation number $\alpha'(c)$.
\QED

\subsection{}\label{mainu}
The asymptotic behaviour of viscosity solution
is described by the following theorem, obtained in \cite{Be}
(see also \cite{BR} for a better proof, and see
\cite{F3} and \cite{Roque} for related results).
Let $T(c)\in \Nm$ be defined by:\\
T(c)=1 if $\rho(c)$ is irrational,\\
T(c)=q if $\rho(c)$ is the rational $p/q$ in lowest terms.\\
\textsc{Theorem }
\begin{itshape}
For each $u\in C(\Tm,\Rm)$,
the sequence
$\big(V^c\big)^{nT(c)}(u)$
is converging uniformly to a fixed point 
of $\big(V^c\big)^{T(c)}$.
\end{itshape}
\section{Entropy solutions and characteristics}
\label{E}

The relation between classical solutions
of Burgers equation and the Hamiltonian dynamics 
is quite well understood from
\ref{char}.
We shall now describe the main properties 
of entropy solutions, with emphasis on their relation
with dynamics.
We will also prove Theorem \ref{main2} and corollary \ref{cor}.

\subsection{}\label{entropy}
A function $y:[t_0,\infty)\lto \Rm$
is called an entropy solution of Burgers equation $(B)$ if :\vs\\
$i.$ The functions $y$ and $H(t,x,y(t,x))$
are locally integrable and 
the equation holds in the sense of distributions :
$$
\int _{[t_0,\infty)\times\Tm }
y(t,x)\partial_t \phi(t,x) +H(t,x,y(t,x))\partial _x\phi(t,x) dtdx
+
\int _{\Tm}y(t_0,x)\phi(t_0,x)dx =0
$$
for all smooth function $\phi:[t_0,\infty)\times\Tm\lto \Rm$
with compact support (see \cite{S} for details).
Note that  the space average 
$
\int _{\Tm}y(t,x)dx
$
is then a constant $c$.\vs\\
$ii.$ The Oleinik inequalities
$$
y_t(x+\delta)-y_t(x)\leq K(t)\delta
$$
hold for all $t>t_0$, $x\in \Tm$
and $\delta >0$,
with a positive and decreasing function $K(t)$.\\
\subsection{}\label{propentropy}
\textsc{Proposition }
\begin{itshape}
For each $t_0\in \Rm$ and each $y_0\in L^{\infty}(\Tm)$,
there exists a unique entropy solution 
$y:[t_0,\infty)\times \Tm\lto\Rm$
of Burgers equation
satisfying
$y(t_0,.)=y_0$.
We have 
$$
y\in L^{\infty}([t_0,\infty)\times\Tm)\cap 
C\big([t_0,\infty),L^1(\Tm)\big).$$
This solution 
is given by $y(t,x)=c+\partial_x u(t,x)$,
where $c=\int y_0$
and $u(t,x)$ is the viscosity solution
of $(HJc)$
of initial condition 
$u_{t_0}(x)=\int_0^x(y_0-c).
$
\end{itshape}
\vs\\
Note that we could associate a function 
$y\in C(]t_0,\infty[,L^1)$ to any initial condition
$y_0\in L^1$ in exactly the same way.
However we do  not have $y\in  C([t_0,\infty[,L^1)$
if $y_0$ is not bounded.
\vs

\begin{proof}
Let us first deal with uniqueness.
The standard method to prove uniqueness
is to use the Oleinik inequalities
\ref{entropy} ii.,
via a duality method, see \cite{H}, Theorem 2.2.1,
or \cite{S}, 2.8.
We shall use the Hamilton-Jacobi equation.
Indeed, let $y(t,x) : [t_0,\infty)\times \Tm\lto \Rm$
be an entropy solution.
Note that this function is locally bounded
in $]t_0,\infty)\times \Tm$
in view of the Oleinik inequalities.
Define  $\tilde u(t,x)=\int_{0}^x(y(t,x)-c)dx-\int_{0}^1(y(t,x)-c)dx,$
where $c=\int y$.
We have, in the sense of distributions,
$\partial_{tx}\tilde u=\partial_t y
=-\partial_x(H(t,x,y))$.
Hence
the distribution
$\partial_t \tilde u(t,x)+H(t,x,y(t,x))$
does not depend on $x$,
and is the locally integrable  function
$f(t)=\int _{\Tm}H(t,x,y(t,x))dx$.
The function $u(t,x)=\tilde u(t,x)-\int_{t_0}^t f(s)ds$,
satisfies
$\partial_t u+ H(t,x,y)=0$
in the sense of distributions,
hence it is locally Lipschitz on $]t_0,\infty)\times \Tm$
and  the Hamilton-Jacobi equation 
$\partial_t u+ H(t,x,c+\partial_x u)=0$
holds almost everywhere.
In addition, the function $u$
satisfies condition $ii.$ of proposition
\ref{Vprop} as a consequence of the Oleinik inequalities
$ii.$ above, 
so that $u$ has to be \textit{the} viscosity solution
of $(HJc)$ of initial condition $u_{t_0}$.
So the only candidate to be an entropy solution
of $(B)$ is $y(t,x)=c+\partial_{x} u(t,x)$.

It is classical to obtain the existence
of entropy solutions as limits
of regular solutions of the viscous equation
$(B_{\mu})$.
However we shall use once more the Hamilton Jacobi equation,
\textit{i.e.} we shall prove that the function
$y(t,x)=c+\partial_x u(t,x)$ introduced in the
discussion on uniqueness is indeed an entropy solution.
Recall that $u_{t_0}$ is Lipschitz hence
$u$ is Lipschitz,
hence $\partial _t u$ and $y=c+\partial _x u$
are well defined in $L^{\infty}([t_0,\infty)\times \Tm)$,
as well as $H(t,x,y(t,x))$.
It is straightforward that $i.$ is satisfied,
and $ii.$ follows
from the property $ii.$ of Proposition \ref{Vprop}.
Since the the injection $L^{\infty}\lto L^1$ is compact,
the continuity of $t\lmto y_t$ follows from the continuity 
of $t\lmto u_t$ and the fact that the functions
$y_t$ are equibounded.
\end{proof}

 \subsection{}\label{Eprop}
We call $E$ the operator $y_0\lmto y_1$,
which can be extended to an operator
$$
E:L^1(\Tm)\lto L^1(\Tm).$$
Let us recall the major properties of $E$.\\
$i.$ The operator $E$ is nondecreasing,
and it is a contraction in $L^1$:
$$
\|E(y)-E(z)\|_1\leq \|y-z\|_1
$$
$ii.$
The operator $E$ is compact. More precisely,
for all $C>0$, there exists $K>0$ such that,
when  $|\int y|\leq C$,
we have $\|E(y)\|_{\infty}\leq K$, and the Oleinik inequality
$$
E(y)(x+\delta)\leq E(y)(x)+K\delta
$$
for all $x$ and $\delta>0$.\\
$iii.$
The operator $E$ preserves the average $c=\int y$.
\vs\\
The contraction property is quoted here for completeness.
It can be obtained as a consequence of the Oleinik
inequalities using a duality method, see \cite{S}.
It could also be proved using the function $u$,
but we shall not write down this proof here.
We shall only use this property in the case where
all functions $y$, $z$, $E(y)$ and $E(z)$ are continuous.
In this very easy case, the contraction property is a direct 
consequence of the area-preservation property of the diffeomorphism
$\phi$, as will become clear in the sequel.

\subsection{}
Consider an entropy solution 
$y(t,x):[t_0,\infty)\times \Tm \lto \Rm$
of Burgers equation,
let $c=\int y$  and let 
$u(t,x):[t_0,\infty)\times \Tm \lto \Rm$
be a viscosity solution
of $(HJc)$ such that $y=c+\partial_x u$.
A trajectory  $\gamma(s):[t_0,\infty)\supset [t,t'] \lto \Tm \times \Rm$
of the Hamiltonian vector-field $X$ is called a 
$y$-characteristic if $\gamma(s)\in \overline{G(y_s)}$ for 
each $s\in [t,t']$.
The following theorem extends the method 
of characteristics to entropy solutions :
\vs\\
\textsc{Theorem (Characteristics) }
\begin{itshape}\\
i. 
Let us fix $t_0\leq t<t'$.
A curve $\gamma(s):]t_0,\infty)\supset[t,t']\lto \Tm\times \Rm$ 
is a $y$-characteristic if and only if
it is calibrated by $u$.\\
ii.
For every $t>t_0$ and every $(x,p)\in \overline{G(y_t)},$
there exists a unique $y$-characteristic 
$\gamma(s):]t_0,t]\lto \Tm\times\Rm$
such that $\gamma(t)=(x,p)$.\\ 
iii.
Let $\gamma(s)=(x(s),p(s)):]t_0,t]\lto \Tm \times \Rm$
be a $y$-characteristic, then
for each $s\in ]t_0,t[$, the function
$y_s$ is continuous at $x(s)$
and $y(s,x(s))=p(s)$.\\
iv.
If in addition $y_{t_0}\in \mE$, then 
for every $t>t_0$ and every $(x,p)\in \overline{G(y_t)},$
there exists a unique $y$-characteristic 
$\gamma(s)=(x(s),p(s)):[t_0,t]\lto \Tm\times\Rm$
such that $\gamma(t)=(x,p)$.
In addition, the function $y_{t_0}$ is continuous
at $x(t_0)$, and $y(t_0,x(t_0))=p(t_0)$.\\ 
\end{itshape}

\subsection{}\label{Edyn}
In terms of the dynamics,
using the notations of \ref{graph}, this theorem implies that
$$
\overline{G(y_{t'})}\subset \phi_{t,t'}
\big(G(y_t)\big)
$$
when $t_0<t<t'$, and that
$$
\overline{G(E(y))}\subset \phi
\big(G(y)\big)
$$
when $y\in \mE$.

\subsection{}\label{Erotation}
Let $y(t,x):[t_0,\infty)\times \Tm \lto \Rm$
be an entropy solution of Burgers equation, 
and let $c=\int y$.
We say that an absolutely continuous curve
$x(s):[t_0,\infty)\supset[t,t']\lto \Tm$
is a weak $y$-characteristic
if it satisfies the equation
$\dot x(s)\in [y^-_s(x(s)),y^+_s(x(s))]$ almost everywhere.
In view of the Oleinik inequalities,
it follows from extended Cauchy-Lipschitz results,
see \cite{H}, Theorem 1.4.2, that
for each $t\geq t_0$ and $x_0\in\Tm$,
there exists one and only one
weak characteristic $x(s):[t,\infty)\lto\Tm$
satisfying  $x(t)=x_0$.
It is clear that $y$-characteristics
are weak $y$-characteristics, hence
Proposition \ref{rotationu} implies:\\
\vs
\textsc{Proposition}
\begin{itshape}
All weak $y$-characteristics  $x(s):[t,\infty)\lto\Tm$
have the same rotation number
$$
\lim \frac{1}{t}\int _{t_0}^t \dot x=\rho(c)=\alpha'(c).
$$
\end{itshape}

\subsection{}
Let us define the operator 
$
U:L^1(\Tm)\lto 
C(\Tm,\Rm)
$
by 
$$
U(y)(x)=\int _0 ^x y(x) dx -x \int _{\Tm}y(x)dx.
$$
It will also be useful to consider 
the operator 
$
\tilde U:
 L^1(\Tm)\lto
\Rm\times C(\Tm)
$
defined by 
$
\tilde U(y)=(\int_{\Tm} y,U(y)).
$
We then have
$$
E=
\tilde U^{-1}\circ
\tilde V\circ
\tilde U.
$$
In order to give a meaning to this expression,
we have to notice that,
if $u$ is in the image of 
$V$,
then $u$ is  Lipschitz
so that $\tilde U^{-1}(c, u)$
is well defined for all $c$ and given by 
$$
\tilde U^{-1}(c,u)(t,x)=
c+\partial_x u(t,x).
$$
In addition, the restriction of $\tilde U^{-1}$
to the image of $\tilde V$
is continuous.
Indeed, let us consider a sequence 
$c_n\lto c$ of reals and a sequence
$u_n$ in the image of $V^{c_n}$,
such that $u_n\lto u$.
The sequence $u_n$ is then equilipshitz,
hence $\partial_x u_n \lto \partial _x u$
in $L^1$.
Theorem \ref{main2} is now a direct consequence of 
Theorem \ref{mainu}.

\subsection{}\label{reversed}
Let us turn our attention to
the dynamics of $\phi$.
It is useful to consider
the reversed Hamiltonian
$\breve H(t,x,p)=H(-t,x,-p)$,
which clearly also satisfies the hypotheses
\ref{Hhyp}.
We denote by $\breve E$ the associated entropy
operator.
Denoting by $S$ the symmetry
$(x,p)\lmto(x,-p)$,
we have for the associated vector-field
$$
\breve  X(t,x,p)
=-S(X(-t,x,-p)).
$$
Hence the flow satisfies
$$
S\circ \breve {\phi} = \phi^{-1}\circ S,
$$
so that 
$$
\overline{G(\breve E(-y))}\subset \phi^{-1}
\big(G(-y)\big)
$$
when $-y\in \mE$.
\vs\\
\textsc{Lemma}
\begin{itshape}
If there exists $n>2$
such that  the function $E^n(y)$ is continuous,
then we have
$$\breve E^{k}(-E^n(y))=-E^{n-k}(y)$$
for $0<k<n$, hence
the functions $E^{n-k}(y)$
are $K$-Lipschitz.
\end{itshape}
\vs\\
\proof
By recurrence, it is enough to prove the Lemma for $k=1$.
We have, noting $y_i=E^i(y)$,
$$
G(-\breve E(-y_n))
\subset 
\phi^{-1}G(y_n)
\subset 
G(y_{n-1})
$$
hence 
$y_{n-1}=-\breve E(-y_n)$.
It follows that $-y_{n-1}$
satisfies the inequalities \ref{Eprop}, ii.
so that $y_{n-1}$ is $K$-Lipschitz since
it is also satisfies these inequalities.
\QED

\subsection{}\label{coroproof}
We are now in a position to prove 
Corollary \ref{cor} for $\phi$.
Let us set $c=\int y$.
The functions $E^n(y)$ are 
$K$-Lipschitz for $n>0$.
In view of Theorem \ref{main2},
the sequence  
$E^{nT(c)}(y)$ has a limit $z$
which is also $K$-Lipschitz
and satisfies $E^{T(c)}(z)=z$.
Let us set for simplicity $y_n=E^{nT(c)}(y)$.
It follows from the Lemma that 
$y_{n-k}=\breve E^{kT(c)}(y_n)$
and $z=\breve E^{T(c)}(z)$.
The contraction property
of $\breve E$,  see \ref{Eprop} $i.$,
implies that the sequence  $\|y_n-z\|_1$
is non decreasing, hence it is identically zero,
so that $y=z$.
\QED

\end{document}